# Autonomous Orbital Correction for Nano Satellites Using J2 Perturbation and LSTM Networks

**Mahya Ramezani[a], Mohammadamin Alandihallaj[b]\* and Andreas M. Hein[c]**

*Interdisciplinary Centre for Security, Reliability and Trust (SnT), University of Luxembourg,*
*2 Avenue de l'Universite L, 4365 Esch-sur-Alzette, Luxembourg;*
*[a] mahya.ramezani@uni.lu (M.R.), [b] amin.hallaj@uni.lu (M.A.), [c] andreas.hein@uni.lu (A.M.H.)*
\* Corresponding Author

**Abstract**

CubeSats offer a cost-effective platform for various space missions, but their limited fuel capacity and susceptibility to environmental disturbances pose significant challenges for precise orbital maneuvering. This paper presents a novel control strategy that integrates a $J_2$-optimized sequence with an LSTM-based low-level control layer to address these issues. The $J_2$-optimized sequence leverages the Earth's oblateness to minimize fuel consumption during orbital corrections, while the LSTM network provides real-time adjustments to compensate for external disturbances and unmodeled dynamics. The LSTM network was trained on a dataset generated from simulated orbital scenarios, including factors such as atmospheric drag, solar radiation pressure, and gravitational perturbations. The proposed system was evaluated through numerical simulations, demonstrating significant improvements in maneuver accuracy and robustness compared to traditional methods. The results show that the combined system efficiently reduces miss distances, even under conditions of high uncertainty. This hybrid approach offers a powerful and adaptive solution for CubeSat missions, balancing fuel efficiency with precise orbital control.

**Keywords:** CubeSat; Orbital dynamics; Orbital maneuvering; LSTM; Uncertainty

**Acronyms/Abbreviations**
**ADS-B** – Automatic Dependent Surveillance-Broadcast
**AI** – Artificial Intelligence
**COTS** – Commercial-Off-The-Shelf
**GRU** – Gated Recurrent Unit
**GPS** – Global Positioning System
**LSTM** – Long Short-Term Memory
**LEO** – Low Earth Orbit
**LVLH** – Local Vertical Local Horizontal
**MPC** – Model Predictive Control
**MSE** – Mean Squared Error
**RNN** – Recurrent Neural Network
**SGP4** – Simplified General Perturbations 4
**STK** – Systems Tool Kit

## 1. Introduction

The rapid advancement of CubeSats has brought about a paradigm shift in the space industry, offering an affordable and efficient platform for conducting a broad spectrum of scientific and technological missions. Since the inaugural CubeSat launch in 2003 [1], these small satellites have garnered widespread adoption due to their compact form factor, lower mass, and shorter development timelines compared to conventional satellites. The appeal of CubeSats lies largely in their utilization of Commercial-Off-The-Shelf (COTS) components, significantly reducing the cost of space missions and making them accessible to academic institutions, small enterprises, and emerging space-faring nations. CubeSats have been successfully employed in diverse applications, including Earth observation [2, 3], communication, air traffic management [4], ionospheric studies [5], active debris removal [6], wildfire monitoring [7], gravity wave detection [8], and even interplanetary exploration [9].

While CubeSats present notable advantages, operational deployment poses several challenges, particularly in achieving precise orbital maneuvers. CubeSats are frequently launched as secondary payloads, meaning their insertion orbits are determined by the requirements of the primary mission rather than their own [10]. This often results in significant orbital deviations that need to be corrected to achieve mission-specific orbital elements. Such corrections are particularly important in missions involving satellite constellations or formations [11, 12], where precise positioning and coordination are essential.

Much of the recent research in CubeSat operations has focused on constellation design and control techniques. For example, Jaffer et al. [13] optimized CubeSat constellations for air traffic monitoring using ADS-B signals, demonstrating the potential of small satellite networks for global traffic management. Similarly, Gu et al. [14] proposed an enhanced 2D coverage map for satellite communication, improving coverage predictions. Bai et al. [15] introduced a multi-stage thrust control strategy for orbital transfers, enhancing maneuver efficiency, while Zhou et al. [16] contributed a parametric formation control method aimed at cooperative satellite missions.

One of the primary challenges CubeSats face in orbital maneuvers is the limitation imposed by their small onboard fuel capacity. Given their compact size, CubeSats can only carry a minimal amount of propellant, making it crucial to implement highly efficient maneuvering strategies to conserve fuel. Furthermore, CubeSat propulsion systems typically generate low thrust, requiring prolonged operational periods to achieve necessary velocity changes (ΔV) [17]. These limitations are compounded by the restricted power availability from CubeSats' small solar arrays, which further limits propulsion system usage. Additionally, advanced control methods, such as Model Predictive Control (MPC) [18-20] are often impractical for CubeSats due to their limited onboard computational resources [21].

To overcome these constraints, various fuel-saving techniques have been proposed in the literature, including exploiting aerodynamic forces, solar radiation pressure, and third-body gravitational effects [22-25]. However, the use of such methods is highly dependent on the specific orbital environment and mission objectives. Furthermore, the inherent uncertainties introduced by COTS components, which often do not meet stringent space-grade standards, add another layer of complexity to achieving precise orbital maneuvers [26].

To overcome these challenges, machine learning presents a promising solution, offering the capability to adapt to unmodeled dynamics and uncertainties [21]. Several studies have demonstrated the effectiveness of machine learning techniques across various applications including path planning for aerial robots [27-29], satellite control [30], and failure prediction [31], showcasing their potential to enhance the accuracy and reliability of complex systems in real-world scenarios. This paper proposes a novel approach that integrates advanced machine learning techniques with traditional orbital mechanics to optimize CubeSat maneuvering. The core of this approach is the development of a machine learning-based maneuvering framework that leverages the Long Short-Term Memory (LSTM) network for real-time trajectory adjustments. LSTMs, known for their ability to manage sequential data and long-term dependencies, are well-suited for the time-series nature of orbital maneuvering data [27, 32]. Unlike previous work that focused on using less computationally intensive models such as GRUs [33, 34], this paper emphasizes the LSTM model's superior ability to capture complex, time-varying orbital dynamics, making it a more effective choice for high-precision orbital corrections.

In the proposed framework, the LSTM model plays a central role in dynamically adjusting the CubeSat's orbital maneuvers. The LSTM continuously monitors discrepancies between the planned trajectory and actual execution, adapting to unmodeled dynamics and disturbances in real-time. This ability to learn from and react to real-time data makes LSTMs particularly effective in environments characterized by uncertainty and non-linearity, such as space missions. Unlike pre-planned maneuver sequences, the LSTM-based system enables the CubeSat to adapt on-the-fly to unexpected conditions, significantly improving maneuver accuracy and reliability.

The use of LSTM networks in this context provides several key advantages. First, LSTMs excel at handling time-series data, making them ideal for sequential decision-making in orbital maneuvers. Second, they offer the capability to process long-term dependencies, which is critical in tracking cumulative effects of small disturbances over time. Finally, the LSTM's adaptive nature allows it to mitigate uncertainties caused by external perturbations, such as those introduced by the variability of COTS components or environmental forces like solar radiation pressure and atmospheric drag.

This paper's main contributions are twofold:

1. Development of an AI-driven Maneuvering Framework: This work introduces a machine learning framework based on LSTM networks to optimize CubeSat orbital maneuvers. The model continuously adapts to real-time deviations and disturbances, enhancing maneuver precision and reducing the need for extensive pre-mission simulations.

2. Integration of Orbital Mechanics with Machine Learning: The paper presents a novel combination of traditional orbital dynamics with cutting-edge machine learning techniques, offering a flexible and adaptive approach to CubeSat control. By leveraging the LSTM model's sequential learning capabilities, the framework ensures accurate and fuel-efficient orbital corrections in real-time.

The remainder of this paper is structured as follows: Section 2 provides an overview of the orbital dynamics and mathematical models governing CubeSat maneuvering. Section 3 details the machine learning framework, including the LSTM architecture and its application in orbital control. Section 4 presents simulation results, demonstrating the effectiveness of the proposed approach in practical scenarios. Finally, Section 5 concludes with a discussion on the implications of integrating AI into CubeSat mission planning and outlines potential future research directions.

## 2. Orbital Dynamics

Precise control of a CubeSat's orbit necessitates an understanding of orbital dynamics, particularly when considering perturbative forces such as the Earth's oblateness, represented by the $J_2$ term. This section introduces the state-space representation of CubeSat orbital dynamics, which serves as the foundation for developing advanced control strategies.

The state of a CubeSat in orbit can be described by its classical orbital elements: semimajor axis ($a$), eccentricity ($e$), inclination ($i$), right ascension of the ascending node ($\Omega$), argument of perigee ($\omega$), and true anomaly ($\theta$). These six elements form the state vector $x = [a, e, i, \Omega, \omega, u]^T$, which evolves over time due to gravitational and external forces.

To model this evolution, we represent the dynamics in state-space form, where the time derivative of the state vector $\dot{x}$ is expressed as a function of both natural orbital motion and control inputs. This is captured by the following first-order differential equation:

$$\dot{x} = A(x) + B(x)u + d \quad (1)$$

Here,

$$A(x) = \begin{bmatrix} 0 \\ 0 \\ 0 \\ -\frac{3R_e^2 J_2 n}{2p^2}\cos i \\ \frac{3R_e^2 J_2 n}{4p^2}(4 - 5\sin^2 i) \\ \frac{h}{r^2} - \frac{3R_e^2 n J_2}{4p^2}\left(\sin^2 i \left(5 - 3\sqrt{1-e^2}\right) + \left(2\sqrt{1-e^2} - 4\right)\right) \end{bmatrix} \quad (2)$$

represents the natural evolution of the orbital elements due to Earth's gravitational field, particularly the perturbation caused by the J$_2$ effect, which induces secular changes in $\Omega$ and $\omega$. The matrix

$$B(x) = \begin{bmatrix} \frac{2a^2}{h}e\sin\theta & \frac{2a^2}{h}\frac{p}{r} & 0 \\ \frac{1}{h}p\sin\theta & \frac{1}{h}((p+r)\cos\theta + re) & 0 \\ 0 & 0 & \frac{r\cos u}{h} \\ 0 & 0 & \frac{r\sin u}{h\sin i} \\ -\frac{p}{he}\cos\theta & \frac{p}{he}\left(1+\frac{r}{p}\right)\sin\theta & -\frac{r\sin u}{h\tan i} \\ 0 & 0 & -\frac{r\sin u}{h\tan i} \end{bmatrix} \quad (3)$$

describes how the control inputs $u = [f_r, f_c, f_n]^T$ (forces in the radial, along-track, and cross-track directions) affect the orbital elements. The term $d$ accounts for external disturbances, including atmospheric drag, solar radiation pressure, and third-body gravitational forces. This compact state-space formulation will be utilized to design control strategies that optimize CubeSat maneuvers, ensuring both fuel efficiency and precision in achieving mission objectives.

## 3. Integrated Control Methodology for CubeSat Orbital Maneuvering

This section presents an integrated control methodology designed to address the challenges of precise orbital maneuvering for CubeSats, taking into account uncertainties and disturbances encountered during operations. The control framework consists of two interconnected layers: a high-level control layer that uses the J$_2$-optimized Sequence to generate waypoints for reaching a terminal orbital state, and a low-level control layer, which employs an LSTM network to design and refine the maneuvers required to reach each waypoint. This approach enhances trajectory accuracy while minimizing fuel consumption and compensating for unmodeled dynamics and external perturbations.

*3.1. Overview of the Two-Level Control Architecture*

The proposed control architecture consists of two layers: a high-level control layer and a low-level control layer. At the core of the high-level control is the J$_2$-optimized Sequence, which calculates an optimal maneuver plan based on the desired final orbital parameters, initial state, and the mission's time constraints. Inputs to the high-level controller include the initial orbital parameters $x_0$, the target orbital parameters $x_f$, and the time window for completing the maneuvers. The output is a maneuver plan specifying the timing and magnitude of each velocity change ($\Delta V$) across several steps, along with the expected orbital parameters after each step.

While the high-level controller provides an optimal sequence for orbital changes, it does not account for real-world disturbances such as unmodeled dynamics, environmental forces, or hardware inaccuracies. These factors can cause deviations between the planned and actual CubeSat behavior. To address this, a low-level control layer is introduced. This layer is designed to monitor the performance of the high-level maneuvers and calculate real-time adjustments to compensate for errors and uncertainties. The low-level controller leverages an LSTM network, which refines the CubeSat's maneuvers, ensuring each waypoint generated by the high-level controller is accurately reached.

The LSTM network is trained offline using simulated data generated from the J$_2$-optimized Sequence and the CubeSat's dynamic model. During this training, the LSTM learns how to adjust the magnitude and timing of the $\Delta V$ maneuvers to compensate for disturbances, such as atmospheric drag or gravitational anomalies, that could affect the CubeSat's trajectory. The inputs to the LSTM include the current state $x_c$ and the target waypoint provided by the high-level controller. The LSTM outputs the necessary adjustments in $\Delta V$ and timing to ensure the CubeSat reaches the waypoint with minimal deviations, even in the presence of real-time uncertainties.

This two-level architecture provides a robust, adaptive system for CubeSat orbital control. The high-level controller ensures fuel-efficient maneuvering by leveraging natural orbital perturbations (e.g., the J$_2$ effect), while the LSTM-based low-level controller provides real-time adaptability to counteract disturbances. A block diagram of the system is presented in Figure 1, which illustrates the interaction between the high-level controller (Green), the low-level LSTM controller (Blue), a Kalman filter (Red), and the CubeSat dynamics (Purple). This diagram shows how the system components work together to execute the optimal maneuver sequence while continuously adapting to real-time conditions.

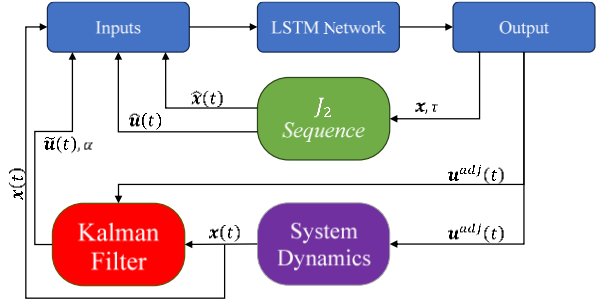

Fig. 1. Schematic diagram of the system.

*3.2. High-Level Control Layer*

The high-level control layer is responsible for generating waypoints and maneuver commands to guide the CubeSat toward its target orbit in a fuel-efficient manner. This layer focuses on leveraging the $J_2$ perturbation to minimize fuel consumption while achieving precise orbital corrections. The maneuvers are divided into in-plane and out-of-plane actions, each aimed at adjusting specific orbital elements.

Orbital Maneuver Planning with $J_2$-Optimized Sequence

To transfer a CubeSat into its desired final orbit, the $J_2$-optimized Sequence is employed, which capitalizes on the natural perturbative effects of the Earth's oblateness ($J_2$). By placing the CubeSat in an appropriate transfer orbit, the $J_2$ perturbation can be used to reduce the errors in the right ascension of the ascending node ($\Omega$) and the argument of latitude ($u$) over the maneuvering period. This method allows the CubeSat to achieve the target orbit with minimal fuel expenditure.

The $J_2$-optimized Sequence works by adjusting the transfer orbit's semi-major axis ($a$) and inclination ($i$) to exploit the $J_2$ perturbation for correcting $\Omega$ and $u$. If the initial and transfer orbits have the same right ascension of ascending node at the start ($\Omega_t = \Omega_0$), the condition for $\Omega$ to converge to the desired value, assisted by the $J_2$ perturbation, is given by:

$$\dot{\Omega}_f = \dot{\Omega}_0 + \frac{\Omega_f - \Omega_0}{\tau} \tag{4}$$

where $\tau$ is the total time of the maneuver sequence. This condition can be substituted into Gauss's variational equation to determine the optimal inclination of the transfer orbit:

$$i_t = \cos^{-1}\left(\frac{-\dot{\Omega}_f + \frac{\Delta\Omega}{\tau}}{3R_e^2 J_2 n_t} 2p_t^2\right) \tag{5}$$

where $R_e$ is the Earth's mean radius, $J_2$ is the second zonal harmonic, $n_t$ is the mean motion of the transfer orbit, and ptp_tpt is the semi-latus rectum of the transfer orbit. The desired semi-major axis $a_t$ of the transfer orbit can be derived by ensuring that $u$ reaches its target value at the end of the maneuver sequence:

$$a_t = a_f \left(1 + \frac{(u_f - u_0) + \delta u_{j2}}{2\pi k}\right)^{\frac{2}{3}} \tag{6}$$

where $k$ represents the number of revolutions of the CubeSat in the transfer orbit, and $\delta u_{j2}$ represents the variation in the argument of latitude due to the $J_2$ perturbation. This variation, $\delta u_{j2}$, is calculated using the following first-order approximation:

$$\delta u_{j2} = -\frac{3R_e^2 J_2}{a_f^3}(\sin^2 i_f - 1)\left(\frac{3\mu}{4a_f^3 n_f} + n_f\right)(a_f - a_t)\tau \\ -\frac{3R_e^2 n J_2}{2a_f^2}\sin(2i_f)(i_f - i_t)\tau \tag{7}$$

This equation captures the secular effects of the $J_2$ perturbation on the CubeSat's argument of latitude during the transfer orbit, accounting for both the semi-major axis difference and the inclination change.

The sequence begins with an in-plane maneuver $\Delta V_c^a(a_t, a_0)$ to adjust the semi-major axis, followed by an out-of-plane maneuver $\Delta \Delta V_n^i(i_t, i_0)$ to modify the inclination. As the CubeSat progresses along the transfer orbit, the $J_2$ perturbation gradually aligns the argument of latitude and the right ascension of the ascending node with the desired values.

Once the argument of latitude reaches the target value, a second in-plane maneuver $\Delta V_c^a(a_f, a_t)$ is executed to adjust the semi-major axis to its final value. Finally, a pure inclination change maneuver $\Delta V_n^i(i_f, i_t)$ is performed at the end of the maneuvering period to bring the CubeSat to the desired inclination.

The total $\Delta V$ required for the $J_2$-optimized Sequence is calculated as:

$$\Delta V_J = |\Delta V_c^a(a_t, a_0)| + |\Delta V_n^i(i_t, i_0)| \\ + |\Delta V_c^a(a_f, a_t)| \\ + |\Delta V_n^i(i_f, i_t)| \tag{8}$$

This sequence allows the CubeSat to reach the desired orbit efficiently by taking advantage of the natural perturbative effects of $J_2$, reducing fuel consumption compared to traditional methods that do not leverage these perturbations.

*3.3. Low-Level Control Layer*

While the high-level control layer generates an optimal maneuver plan using the $J_2$-optimized Sequence, real-world conditions like external disturbances and unmodeled dynamics can cause deviations from the expected trajectory. To address these issues, an LSTM-based low-level control layer is introduced. The LSTM network provides real-time adaptive adjustments to ensure the CubeSat reaches its target orbit accurately despite environmental uncertainties. The LSTM's ability to capture long-term dependencies in sequential data makes it ideal for adjusting maneuvers based on both current and past states.

The LSTM network takes several inputs: the current orbital state $x_c$, the next target waypoint $x_i$ from the high-level control, pointing accuracy estimated via a Kalman filter, and historical state data—a sequence of past states and control inputs. These inputs allow the LSTM to predict necessary adjustments to the velocity $\Delta V$ and the timing of maneuvers.

At each time step $t$, the LSTM updates its internal state based on three key components:

The forget gate $f(t)$ decides how much of the previous cell state $c(t-1)$ to retain:
$$f(t) = \sigma(W_f x(t) + U_f h(t-1) + b_f) \quad (9)$$

The input gate $i(t)$ controls how much new information to add to the cell state, and the candidate cell state $\tilde{c}(t)$ is computed:
$$i(t) = \sigma(W_i x(t) + U_i h(t-1) + b_i)$$
$$\tilde{c}(t) = \tanh(W_c x(t) + U_c h(t-1) + b_c) \quad (10)$$

The cell state $c(t)$ is updated by combining the previous cell state and the new candidate information:
$$c(t) = f(t).c(t-1) + i(t).\tilde{c}(t) \quad (11)$$

The output gate $o(t)$ controls how much of the updated cell state to output, generating the hidden state $h(t)$:
$$o(t) = \sigma(W_o x(t) + U_o h(t-1) + b_o)$$
$$h(t) = o(t).\tanh(c(t)) \quad (12)$$

The final output of the LSTM includes the adjusted $\Delta V$ and timing for the maneuver:
$$y(t) = W_y h(t) + b_y \quad (13)$$

The LSTM is trained offline using simulated CubeSat dynamics under various conditions, including disturbances like atmospheric drag and gravitational perturbations. The objective is to minimize the Mean Squared Error (MSE) between the LSTM's predicted adjustments and the correct adjustments derived from simulations. The loss function is given by:
$$\mathcal{L} = \sum_{i \in N} \left(\lambda |y_p(t_i) - y_t(t_i)|^2\right) \quad (14)$$
where $y_p(t_i)$ is the predicted adjustment, $y_t(t_i)$ is the target adjustment, and $N$ is the number of training samples.

Once trained, the LSTM is deployed onboard the CubeSat, providing real-time adjustments to the high-level maneuver plan. The LSTM ensures the CubeSat remains on the desired trajectory, compensating for disturbances and ensuring precise control throughout the mission.

**4. Experiment Results and Discussion**

This section evaluates the performance of the proposed control strategy through numerical simulations, focusing on the effectiveness of the LSTM-based low-level control layer in adjusting orbital maneuvers under various conditions. The LSTM network was trained using a comprehensive dataset generated in a simulated environment, and its performance is compared against a GRU-based approach.

The LSTM-based low-level control layer was evaluated by training the network on a dataset generated using Systems Tool Kit (STK) 12.3 and MATLAB 2024a. The simulated environment modeled various disturbances, such as atmospheric drag, solar radiation pressure, and gravitational perturbations from the Moon and Sun. These scenarios covered a wide range of possible conditions, providing a robust foundation for training the LSTM network.

To ensure generalization and avoid overfitting, a five-fold cross-validation approach was employed during the training process. This method divides the dataset into five subsets, using four subsets for training and one for validation in each iteration. Through this process, the model was evaluated on unseen data, ensuring reliable performance in real-world conditions.

The optimal hyperparameters for the LSTM network were selected based on multiple simulation trials. The final configuration is outlined in Table 1.

Table 1. Optimal hyperparameters for the GRU network.

| Hyperparameter | Value |
| --- | --- |
| Number of Layers | 3 |
| Number of Neurons | 128 per layer |
| Learning Rate | 0.001 |
| Batch Size | 64 |
| Dropout Rate | 0.4 |
| Gradient Clipping | 0.5 |
| Activation Function | ReLU |
| Optimizer | Adam |

The LSTM network was trained over multiple iterations, with the loss function, defined as Mean Squared Error (MSE).

In addition to evaluating the accuracy of the LSTM-based approach, its computational efficiency was compared to a GRU-based method. The total time required for training and prediction was recorded for both networks. As shown in Table 2, the LSTM network required more computational resources compared to the GRU network, but this was offset by its superior handling of long-term dependencies in the orbital maneuvering data.

Table 2. Average of computational time comparison between GRU and LSTM networks for one epoch.

| Method | Training Time (s) | Prediction Time (ms) |
| --- | --- | --- |
| GRU | 743 | 11 |
| LSTM | 995 | 18 |

The LSTM network's slightly longer training and prediction times are attributable to its more complex architecture. However, the LSTM offers better accuracy for tasks requiring the handling of longer sequential dependencies, making it more suitable for CubeSat control applications.

To assess the robustness of the LSTM-based control system, a Monte Carlo simulation was conducted. The initial conditions for the CubeSat's altitude were varied randomly between 700 km and 800 km, and the pointing accuracy was modeled as a normal distribution with a covariance of $\alpha$ degrees, simulating uncertainties in the thrust direction.

The accuracy of the system was evaluated by measuring the miss distance, defined as the minimum distance between the CubeSat and the target orbit during the maneuvering period, as a function of $\alpha$. Figure 2 presents the miss distance for the system using only the $J_2$-optimized sequence without the LSTM-based control. As expected, the miss distance increases significantly with higher values of $\alpha$, highlighting the limitations of relying solely on pre-planned sequences under conditions of uncertainty.

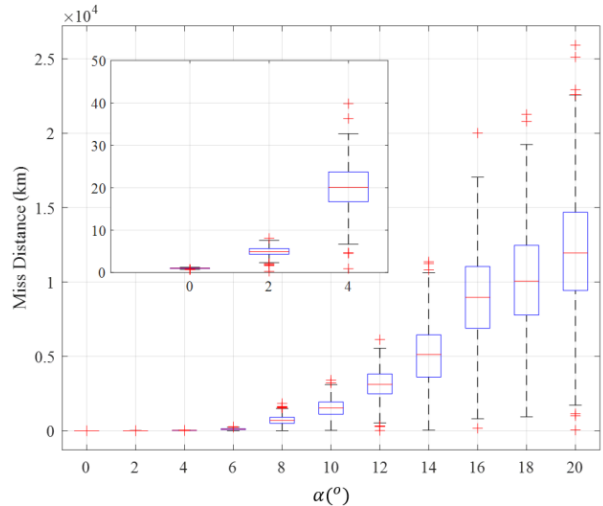

Fig. 2. Miss distance as a function of pointing accuracy $\alpha$ for the system using only the $J_2$-optimized sequence.

In contrast, Figure 3 illustrates the performance of the proposed combined system, which integrates the $J_2$-optimized sequence with the LSTM network and Kalman filter. The results clearly show that the LSTM-based controller effectively compensates for thrust direction uncertainties, maintaining a much lower miss distance even at higher values of $\alpha$. For instance, with $\alpha = 20$ degrees, the average miss distance remains around 10 km, a significant improvement over the system without the LSTM.

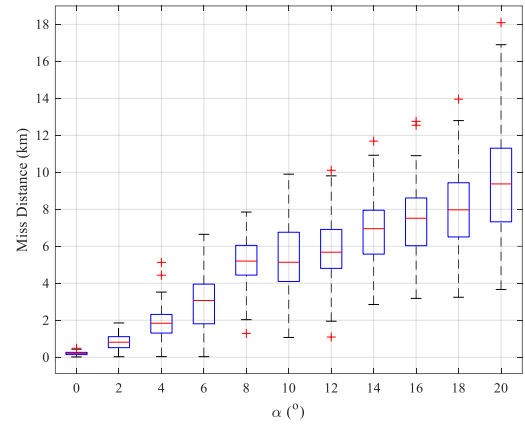

Fig. 3. Miss distance as a function of pointing accuracy $\alpha$ for the combined system using the $J_2$-optimized sequence, LSTM, and Kalman filter.

The results from the Monte Carlo simulation demonstrate that the LSTM-based control system is highly effective in handling real-time uncertainties. By adapting to real-time data and leveraging the predictive power of the LSTM network, the system is able to maintain trajectory accuracy and robustness, even in dynamic and uncertain environments.

The results underscore the benefits of combining the $J_2$-optimized sequence with machine learning techniques like LSTM. While the $J_2$-optimized sequence provides fuel-efficient maneuvers by exploiting natural orbital perturbations, the LSTM enhances the system's ability to adapt to real-world conditions. The Monte Carlo simulation shows that the LSTM-based control system significantly improves performance, particularly in the presence of thrust direction uncertainties. This hybrid approach effectively balances fuel efficiency with robust control, making it highly suitable for CubeSat missions where precision and adaptability are critical.

## 5. Conclusions

This paper presented a novel approach to CubeSat orbital maneuvering by integrating a $J_2$-optimized sequence with an advanced LSTM-based low-level control layer. The proposed control architecture was designed to handle real-time uncertainties and external perturbations that can affect the CubeSat's trajectory, while simultaneously optimizing fuel consumption.

The $J_2$-optimized sequence efficiently exploited the Earth's oblateness to reduce the number and magnitude of required maneuvers, resulting in significant fuel savings. However, to account for real-world uncertainties such as environmental disturbances and hardware imperfections, the LSTM network was introduced as part of the low-level control layer. This network dynamically adjusted the velocity and timing of maneuvers, adapting to deviations from the planned trajectory in real time.

Through a comprehensive training process, using a dataset generated from various orbital scenarios, the LSTM demonstrated its ability to accurately predict the necessary corrections.

Simulation results confirmed the effectiveness of the combined system, particularly under conditions of high uncertainty. The LSTM-based control layer, supported by a Kalman filter for state estimation, significantly reduced miss distances, even in the presence of substantial thrust misalignments and pointing inaccuracies. In comparison to traditional methods, the LSTM network provided greater robustness, accuracy, and adaptability, making it well-suited for CubeSat missions operating in dynamic and unpredictable environments.

The integration of the $J_2$-optimized sequence and machine learning represents a promising step toward more efficient and reliable orbital control strategies for small satellites. This hybrid approach not only improves maneuver accuracy but also enhances mission success by ensuring that CubeSats remain on course despite environmental uncertainties. Future work could explore further optimizations of the LSTM network, including the integration of additional environmental factors and the potential for onboard real-time learning capabilities to further increase system resilience.

*List of references*

[1] M. Swartwout, "The first one hundred cubesats: A statistical look," *Journal of small Satellites,* vol. 2, no. 2, pp. 213-233, 2013.

[2] S. Nag, A. S. Li, and J. H. Merrick, "Scheduling algorithms for rapid imaging using agile Cubesat constellations," *Advances in Space Research,* vol. 61, no. 3, pp. 891-913, 2018.

[3] M. A. Alandihallaj and M. R. Emami, "Multiple-payload fractionated spacecraft for earth observation," *Acta Astronautica,* vol. 191, pp. 451-471, 2022.

[4] S. Wu, W. Chen, C. Cao, C. Zhang, and Z. Mu, "A multiple-CubeSat constellation for integrated earth observation and marine/air traffic monitoring," *Advances in Space Research,* vol. 67, no. 11, pp. 3712-3724, 2021.

[5] Y. Duann *et al.*, "IDEASSat: A 3U CubeSat mission for ionospheric science," *Advances in Space Research,* vol. 66, no. 1, pp. 116-134, 2020.

[6] H. Hakima, M. C. Bazzocchi, and M. R. Emami, "A deorbiter CubeSat for active orbital debris removal," *Advances in Space Research,* vol. 61, no. 9, pp. 2377-2392, 2018.

[7] M. Alandihallaj and M. R. Emami, "Monitoring and Early Detection of Wildfires Using Multiple-payload Fractionated Spacecraft," in *73rd International Astronautical Congress (IAC)*, 2022.

[8] J. Westerhoff *et al.*, "LAICE CubeSat mission for gravity wave studies," *Advances in Space Research,* vol. 56, no. 7, pp. 1413-1427, 2015.

[9] T. Kohout *et al.*, "Feasibility of asteroid exploration using CubeSats—ASPECT case study," *Advances in Space Research,* vol. 62, no. 8, pp. 2239-2244, 2018.

[10] T. Villela, C. A. Costa, A. M. Brandão, F. T. Bueno, and R. Leonardi, "Towards the thousandth CubeSat: A statistical overview," *International Journal of Aerospace Engineering,* vol. 2019, no. 1, p. 5063145, 2019.

[11] M. A. Alandihallaj and M. R. Emami, "Satellite replacement and task reallocation for multiple-payload fractionated Earth observation mission," *Acta Astronautica,* vol. 196, pp. 157-175, 2022.

[12] M. R. Emami and M. A. Alandihallaj, "Performance Enhancement of Fractionated Spacecraft for Earth Observation," *44th COSPAR Scientific Assembly. Held 16-24 July,* vol. 44, p. 57, 2022.

[13] G. Jaffer *et al.*, "Air traffic monitoring using optimized ADS-B CubeSat constellation," *Astrodynamics,* vol. 8, no. 1, pp. 189-208, 2024.

[14] Y. Gu, Y. Chen, Y. Zhang, G. Wu, and S. Bai, "Extended 2D Map for Satellite Coverage Analysis Considering Elevation-Angle Constraint," *IEEE Transactions on Aerospace and Electronic Systems,* 2024.

[15] S. Bai, Y. Wang, H. Liu, and X. Sun, "Finite-Thrust Lambert Transfer Based on Multistage Constant-Vector Thrust Control," *IEEE Transactions On Aerospace and Electronic Systems,* vol. 59, no. 5, pp. 4947-4967, 2023.

[16] H. Zhou, B. Jiao, Z. Dang, and J. Yuan, "Parametric formation control of multiple nanosatellites for cooperative observation of China Space Station," *Astrodynamics,* vol. 8, no. 1, pp. 77-95, 2024.

[17] B. Yost and S. Weston, "State-of-the-art small spacecraft technology," 2024.

[18] M. Amin Alandihallaj, N. Assadian, and K. Khorasani, "Stochastic model predictive control-based countermeasure methodology for satellites against indirect kinetic cyber-attacks," *International Journal of Control,* vol. 96, no. 7, pp. 1895-1908, 2023.

[19] M. A. Alandihallaj, N. Assadian, and R. Varatharajoo, "Finite-time asteroid hovering via multiple-overlapping-horizon multiple-model predictive control," *Advances in Space Research,* vol. 71, no. 1, pp. 645-653, 2023.

[20] M. Alandihallaj, B. C. Yalcin, M. Ramezani, M. A. Olivares Mendez, J. Thoemel, and A. Hein, "Mitigating fuel sloshing disturbance in on-orbit satellite refueling: an experimental study," in *International Astronautical Congress IAC*, 2023.

[21] M. Ramezani, M. Atashgah, M. Alandihallaj, and A. Hein, "Reinforcement Learning for Planning and Task Coordination in a Swarm of CubeSats: Overcoming Processor Limitation Challenges," in *International Astronautical Congress*, Baku, Azerbaijan, October 2-4 2023.

[22] X. Shao, M. Song, D. Zhang, and R. Sun, "Satellite rendezvous using differential aerodynamic forces under J2 perturbation," *Aircraft Engineering and Aerospace Technology: An International Journal,* vol. 87, no. 5, pp. 427-436, 2015.

[23] B. Villac and D. Scheeres, "New class of optimal plane change maneuvers," *Journal of guidance, control, and dynamics,* vol. 26, no. 5, pp. 750-757, 2003.

[24] D. Mishne and E. Edlerman, "Collision-avoidance maneuver of satellites using drag and solar radiation pressure," *Journal of Guidance, Control, and Dynamics,* vol. 40, no. 5, pp. 1191-1205, 2017.

[25] D. Mishne, "Collision Avoidance Maneuver of Propulsionless Satellite, Using Solar Radiation Pressure," in *The 56th Israel Annual Conference on Aerospace Sciences, Israel Soc. of Aeronautics and Astronautics Paper WeL2T3*, vol. 3.

[26] A. Sethi, V. Thakurta, N. Gajanur, B. S. Cheela, K. S. Sadasivan, and R. Hosangadi, "Implementation of COTS components for CubeSat applications," in *2017 IEEE Aerospace Conference*, 2017: IEEE, pp. 1-11.

[27] M. Ramezani and M. Amiri Atashgah, "Energy-Aware Hierarchical Reinforcement Learning Based on the Predictive Energy Consumption Algorithm for Search and Rescue Aerial Robots in Unknown Environments," *Drones,* vol. 8, no. 7, p. 283, 2024.

[28] M. Ramezani, M. Atashgah, J. L. Sanchez-Lopez, and H. Voos, "Human-Centric Aware UAV Trajectory Planning in Search and Rescue Missions Employing Multi-Objective Reinforcement Learning with AHP and Similarity-Based Experience Replay," in *2024 International Conference on Unmanned Aircraft Systems (ICUAS)*, 2024: IEEE, pp. 177-184.

[29] M. Ramezani, M. A. A. Atashgah, and A. Rezae, "A Fault-Tolerant Multi-Agent Reinforcement Learning Framework for UAVs-UGV Coverage Path Planning," *Drones,* vol. 8, 2024.

[30] M. Ramezani, M. A. Alandihallaj, and A. M. Hein, "PPO-Based Dynamic Control of Uncertain Floating Platforms in Zero-G Environment," in *2024 IEEE International Conference on Robotics and Automation (ICRA)*, 2024: IEEE, pp. 11730-11736.

[31] M. A. Alandihallaj, M. Ramezani, and A. M. Hein, "MBSE-Enhanced LSTM Framework for Satellite System Reliability and Failure Prediction," in *MODELSWARD*, 2024, pp. 349-356.

[32] M. Ramezani, H. Habibi, J. L. Sanchez-Lopez, and H. Voos, "UAV path planning employing MPC-reinforcement learning method considering collision avoidance," in *2023 International Conference on Unmanned Aircraft Systems (ICUAS)*, 2023: IEEE, pp. 507-514.

[33] R. Dey and F. M. Salem, "Gate-variants of gated recurrent unit (GRU) neural networks," in *2017 IEEE 60th international midwest symposium on circuits and systems (MWSCAS)*, 2017: IEEE, pp. 1597-1600.

[34] M. Ramezani, M. Alandihallaj, and A. M. Hein, "Fuel-Efficient and Fault-Tolerant CubeSat Orbit Correction via Machine Learning-Based Adaptive Control," *Aerospace,* vol. 11, no. 10, p. 807, 2024.